\documentclass[letterpaper,11pt]{article}
\usepackage{amsmath, amsthm, amssymb}    
\usepackage{bridges}			
\usepackage{graphicx}			
\usepackage[colorlinks=true, urlcolor=blue, citecolor=black, linkcolor=black]{hyperref}  
\usepackage[]{algorithm2e}      
\usepackage{subfloat}
\usepackage{subcaption}			
\usepackage{lipsum}

\usepackage[utf8]{inputenc}
\usepackage{amsfonts}
\usepackage{mathtools}
\usepackage{wrapfig}
\usepackage{tikz}
\usepackage{comment}

\newtheorem{theorem}{Theorem}[section]

\newtheorem{lemma}[theorem]{Lemma}
\newtheorem{conjecture}[theorem]{Conjecture}
\theoremstyle{definition}

\theoremstyle{remark}

\title{Eden model for Pentagons}  

\author{\'Erika Rold\'an\textsuperscript{1,2}, Claudia Silva\textsuperscript{1} and Rosemberg Toala-Enriquez\textsuperscript{3}
\vspace{10pt}\\
\textsuperscript{1} ScaDS.AI, Leipzig University; roldan@mis.mpg.de, callame@ciencias.unam.mx\\
\textsuperscript{2} Max Planck Institute for Mathematics in the Sciences; roldan@mis.mpg.de\\
\textsuperscript{3} George Mason University; rtoalaen@gmu.edu
} 

\date{}					

\begin{document}

\maketitle

\thispagestyle{empty}

\begin{abstract}
We study topological and geometric properties of a cell growth process in the Euclidean plane, where the cells are regular pentagons. To explore the aesthetic aspects of this model, we employ a laser cutter on various materials to create physical representations for some simulations of the model.  
\end{abstract}

\section{The Model}
\label{background}
The model consists of gluing regular pentagons along their edges randomly. We start with a regular pentagon centered at the origin. At each stage, a free edge is chosen uniformly (from those available at that time), a new pentagon is attached to it and the list of free edges is updated. Here, free edge means that a pentagon can be glued to it without overlapping with the interior of an already present pentagon. This random process is inspired by the Eden Model which is a First Passage Percolation Process (FPP) based on the regular square tessellation of the plane \cite{eden1961two}. Recently \cite{manin2023topology, hua2023local, damron2024number}, people have studied the topology of FPP processes based on lattices. Here, in the pentagon model, we do not have an underlying lattice and there will be holes that become impossible to cover at any future time. One of the consequences of this fact is that in the Pentagon model, there is also no notion of a convex limiting shape. Moreover, we can observe from our simulations that the perimeter and the number of holes grow linearly with respect to the number of tiles and this shows a clearly different behavior that differs from the classical FPP models based on lattices.

\begin{figure}[h]
\centering
\subfloat[Cell growth model with 200 pentagons]{
        \includegraphics[scale=0.5]{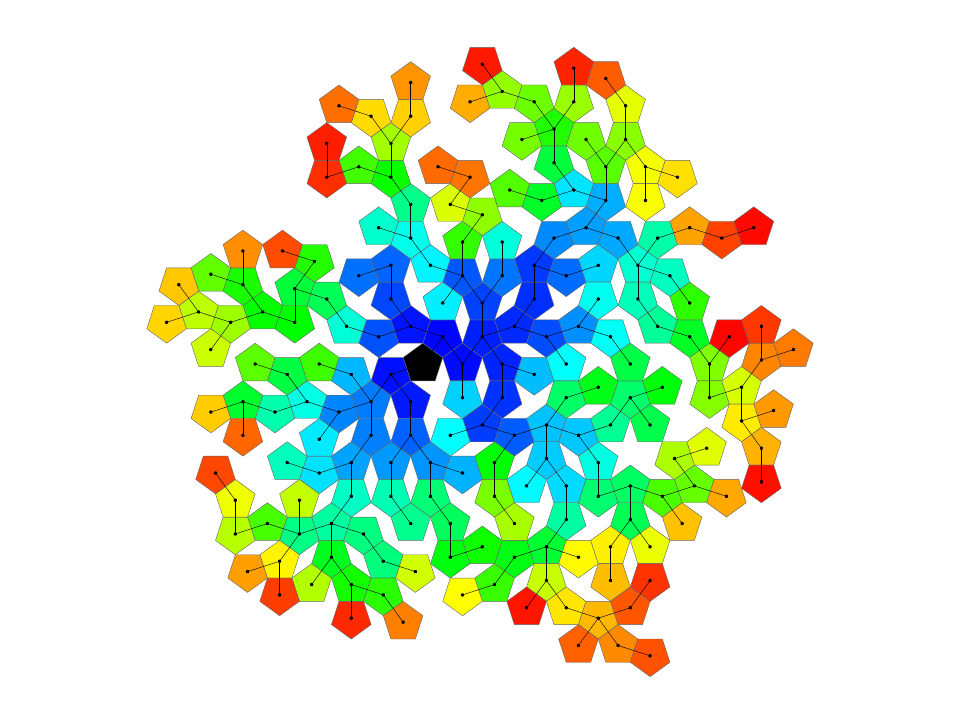}
    }
\subfloat[Model with 10000 pentagons]{
        \includegraphics[scale=0.5]{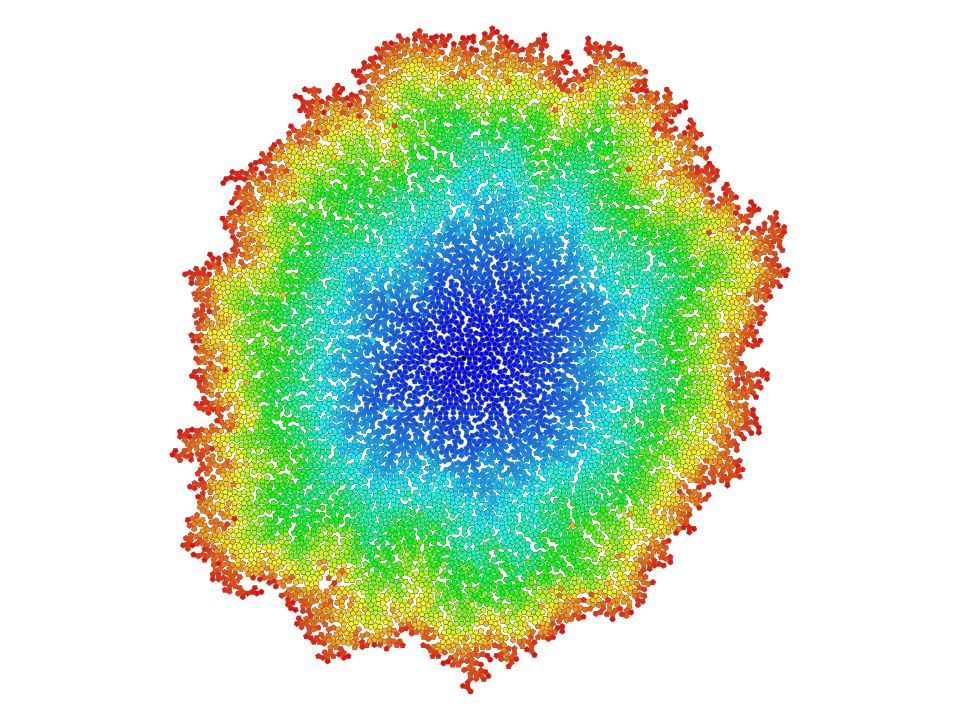}
    } 
\caption{The black pentagon is centered at the origin. The combination of rainbow hues indicates the stage at which that pentagon was placed. The underlying tree represents how the pentagons were attached. }  
\label{pentagons200}
\end{figure}

\section{Geometrical and Topological Properties}

A simple non-intuitive fact is that pentagons come only in two orientations, the first one and its reflection. 

\begin{lemma}
Let $P_0$ be the original pentagon and $P_1$ a reflection of $P_0$ over one side. Then, all pentagons in the model are translations of $P_0$ or $P_1$. Moreover, the sides of all pentagons intersect at angles multiples of $36^{\circ}$.
\end{lemma}
\begin{proof}
   The sides of the first regular pentagon are parallel to the lines of a pencil of 5 lines intersecting at angles of $36^{\circ}$. When a new pentagon is glued, it is obtained by reflecting over the gluing side. The angles do not change under reflection, and one side remains fixed. Therefore, the sides of the new pentagons are still parallel to the same pencil of 5 lines, but the pentagon is now upside down. We conclude that all the pentagons in the model have sides parallel to this pencil of 5 lines and therefore can only be in two possible orientations. 
\end{proof}

\begin{wrapfigure}[13]{r}{0.45\textwidth}
\centering
    \includegraphics{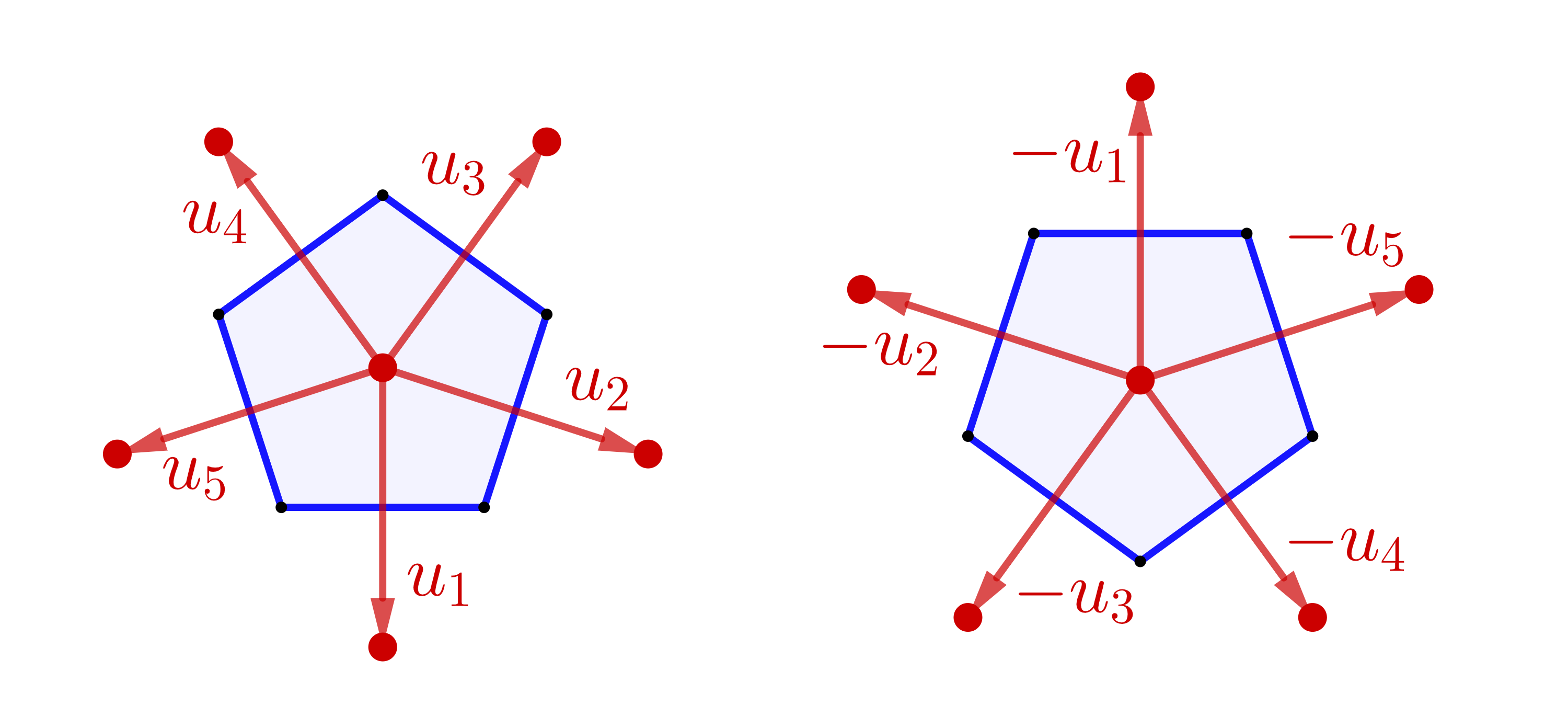}
\caption{The two orientations of the pentagons and the directional vectors of possible neighboring centers.}
\label{fig:orientations}
\end{wrapfigure}

Now we pay attention to the locus of centers of polygons in the model. We notice that the center of a new possible pentagon can be obtained by translating the center of the old pentagon along 5 vectors, see Figure \ref{fig:orientations}. We denote these vectors by $u_1, u_2, u_3, u_4, u_5$. Let us regard $\{u_1, u_2\}$ as a basis for $\mathbb{R}^2$. Then $u_3 = - u_1 + \phi u_2, u_4 = -\phi u_1 - \phi u_2 ,$ and $ u_5 = \phi u_1 - u_2,$ where $\phi = 2 \cos(2\pi/5) = \frac{\sqrt{5}-1}{2}$. Since $\phi$ is irrational, the usual argument shows that the set of integer linear combinations of the vectors $u_1, u_2, u_3, u_4, u_5$ form a dense subset of $\mathbb{R}^2$. This, however, does not take into account the non-overlapping geometric constraint of the pentagons. We conjecture that this difficulty can be overcome. 
\begin{conjecture}
The set of possible centers of pentagons in the model is dense outside a disk around the origin. 
\end{conjecture}

\subsection*{The Graph}
In what follows, we aim to examine the graph formed by the vertices and edges of the pentagons. Our analysis will focus on the growth of the number of vertices, edges, and holes. The \emph{graph} of the model comprises the vertices and edges of the pentagons. If a vertex lies on the edge of another pentagon, the latter is treated as two distinct edges. A \emph{hole} is a bounded component of the complement of the union of all pentagons.   

\begin{figure}[h]
    \centering
\subfloat[Edge-edge]{
        \includegraphics[height=0.7in]{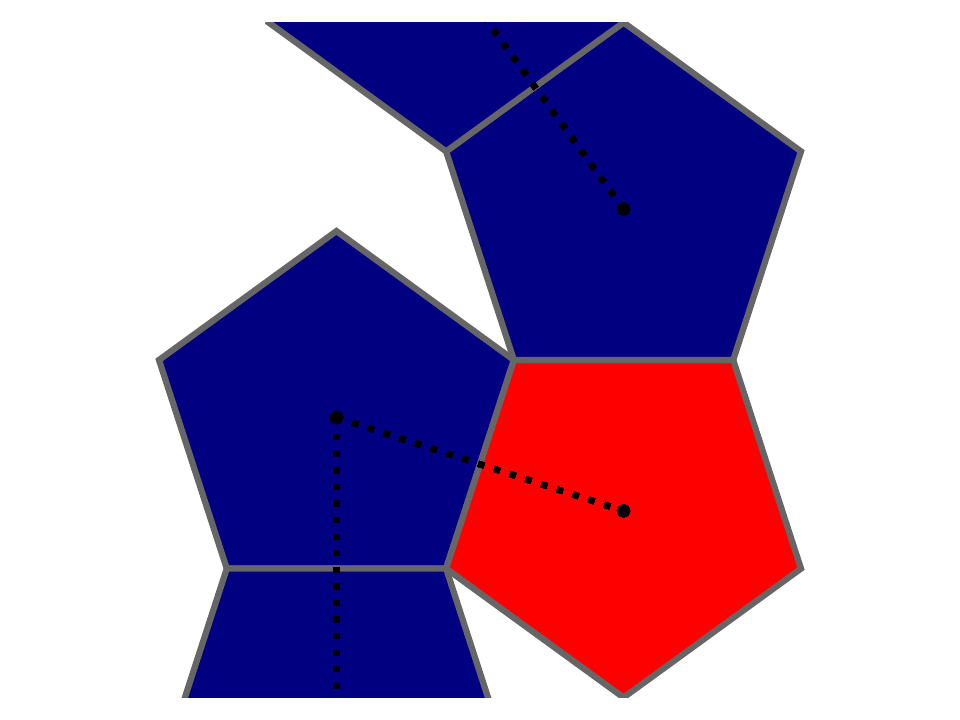}
    }
\subfloat[Partial edge]{
        \includegraphics[height=0.7in]{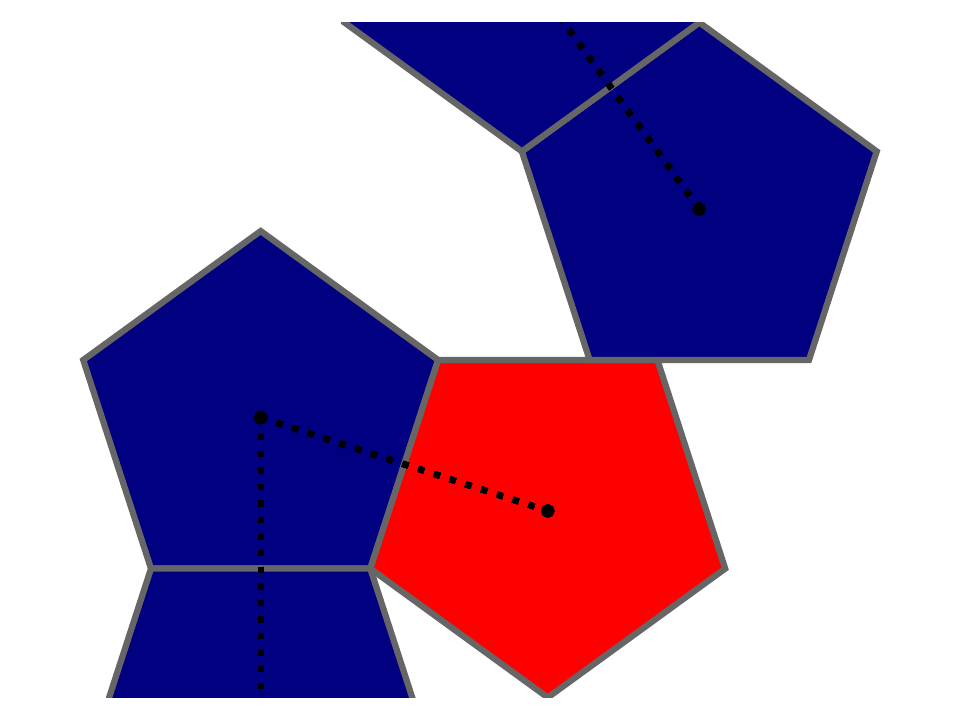}
    }
\subfloat[Vertex-Edge]{
        \includegraphics[height=0.7in]{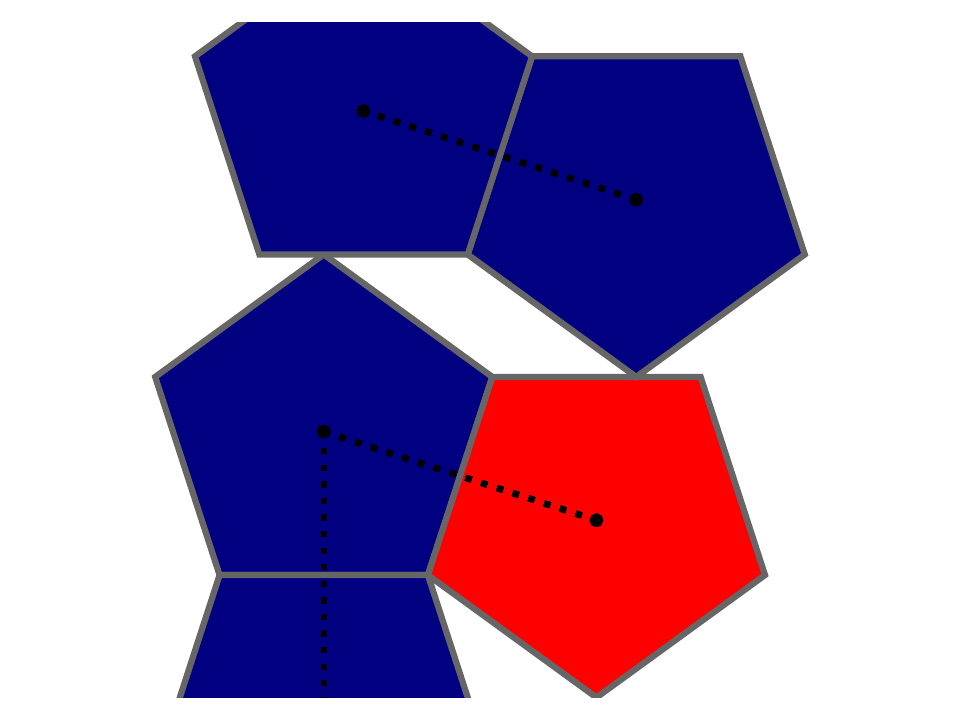}
    }
\subfloat[Vertex-vertex]{
        \includegraphics[height=0.7in]{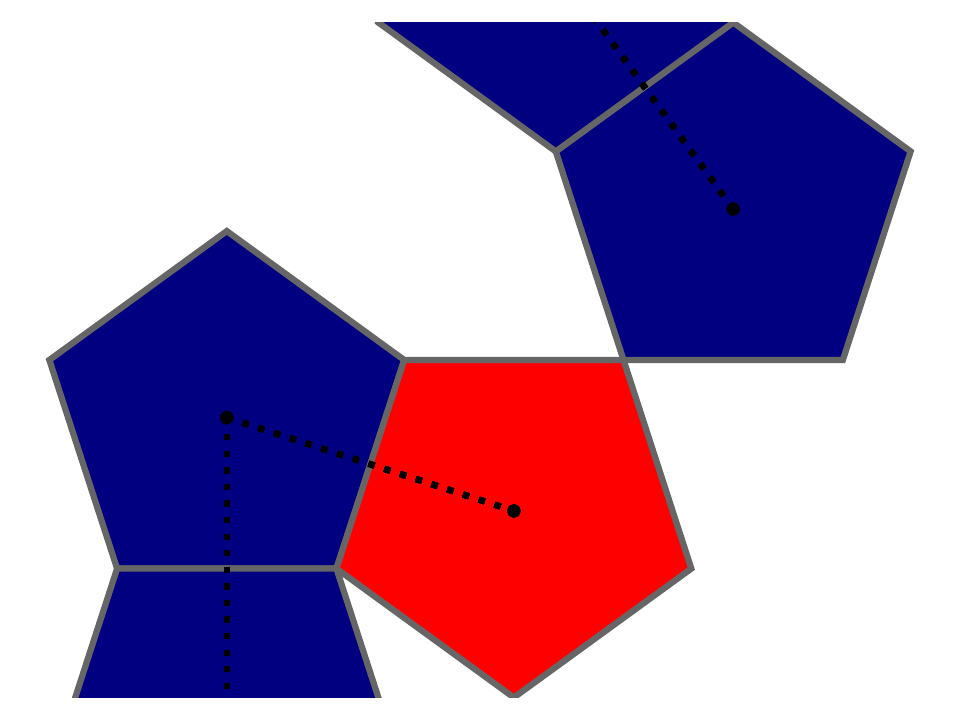}
    }
\caption{A new pentagon (red) can touch the rest of the structure in a combination of 4 different ways.}
\label{fig:touches}
\end{figure}

\begin{figure}[htb]
    \centering
    \includegraphics[width=\textwidth]{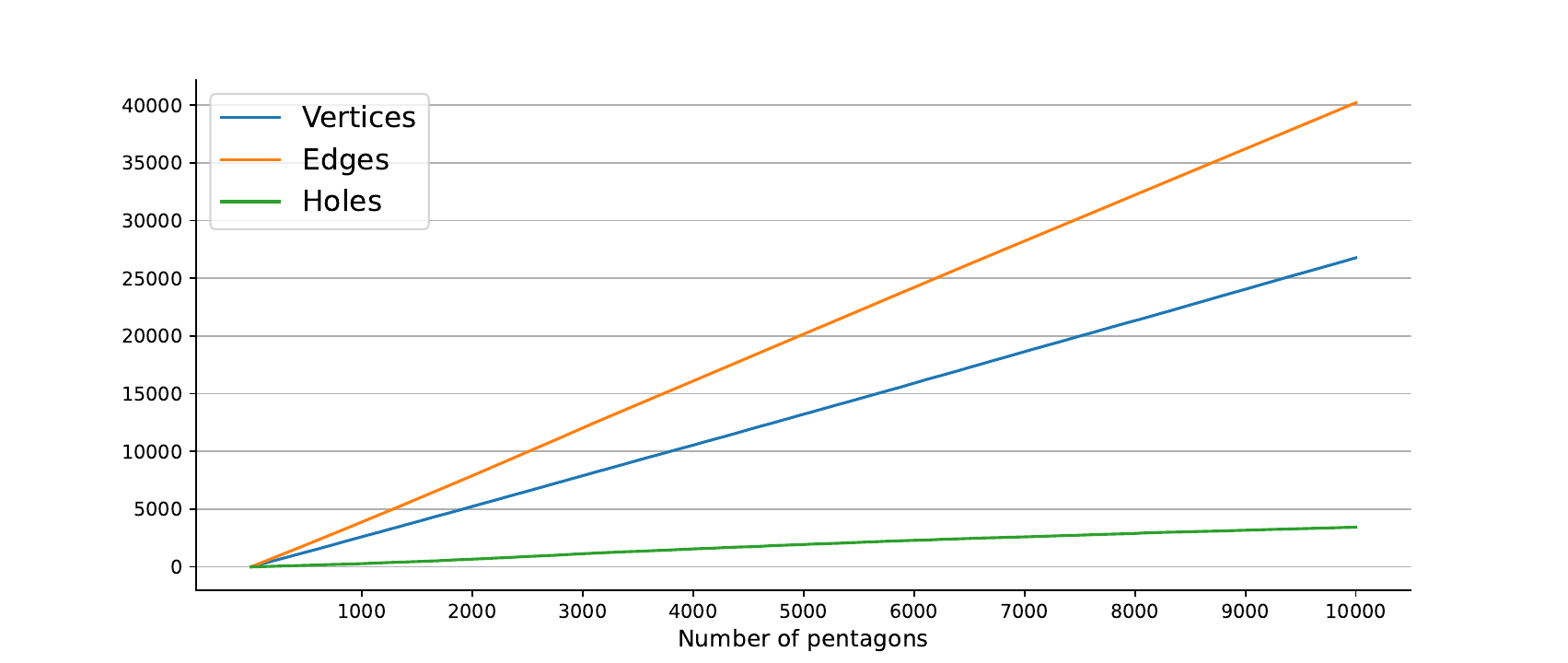}
\caption{Growth of parameters concerning the number of pentagons $n$.}
\label{fig:parameters}
\end{figure}

When a new pentagon is attached, besides the gluing edge, it can touch the existing structure in a combination of 4 different configurations, as shown in Figure \ref{fig:touches}. In each case, when attaching a new pentagon, the increment in the number of vertices, edges, and holes is bounded by 5. We conjecture the following limits, supported by numerical simulations as shown in Figure \ref{fig:parameters}.

\begin{conjecture}
Let $n$ be the number of pentagons in the model, $V(n)$, $E(n)$, and $H(n)$ the number of vertices, edges, and holes, respectively. Then, the expected value of these parameters satisfies the following limits:
\begin{align*}
    \lim_{n\to \infty} \frac{\mathbb{E}[V(n)]}{n} \approx 2.68, \qquad  \lim_{n\to \infty} \frac{\mathbb{E}[E(n)]}{n} \approx 4.02,     \qquad \lim_{n\to \infty} \frac{\mathbb{E}[H(n)]}{n} \approx 0.34.  
\end{align*}
\end{conjecture}

Euler's formula in this context is $V-E+n+H=1$.

\subsection*{Description and Analysis of Holes}
\label{holes}
    
Recall that all the sides of the pentagons form angles multiples of $36^{\circ}$ between themselves. This observation has the consequence that the holes left in the structure are polygons whose angles are multiples of $36^{\circ}$. 

Consider a hole with $l$ sides whose angles are $36^{\circ} a_i$ for some positive integers $a_i$, $i=1,\dots l$. Then, the angle sum condition of the hole polygon is equivalent to: 
    \begin{align}
        a_1 +a_2 + \dots +a_l = 5(l-2).
    \end{align}

\begin{figure}[h!]
    \centering
\subfloat[Triangle]{
        \includegraphics[height=0.87in]{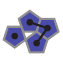}
    }
\subfloat[Arrow]{
        \includegraphics[height=0.87in]{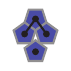}
    }
\subfloat[Paralelogram]{
        \includegraphics[width=0.87in]{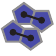}
    }
\subfloat[Diamond]{
        \includegraphics[width=0.87in]{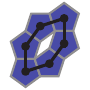}
    }
\subfloat[Ship]{
        \includegraphics[height=0.87in]{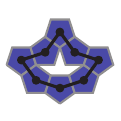}
    }
\subfloat[Three peaks]{
        \includegraphics[height=0.87in]{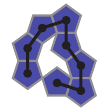}
    }
    
\subfloat[Crown]{
        \includegraphics[height=0.8in]{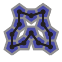}
    }
\subfloat[Double ship]{
        \includegraphics[width=0.8in]{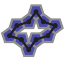}
    }
\subfloat[Pigeon]{
        \includegraphics[height=0.8in]{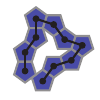}
    }
\subfloat[Triple ship]{
        \includegraphics[height=0.8in]{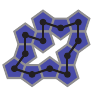}
    }
\subfloat[Claw]{
        \includegraphics[height=0.8in]{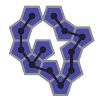}
    }
\subfloat[Whale]{
        \includegraphics[height=0.6in]{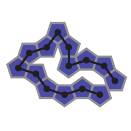}
    }
    
\subfloat[Fox]{
        \includegraphics[height=0.6in]{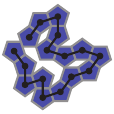}
    }
\subfloat[Bird]{
        \includegraphics[width=0.6in]{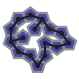}
    }
\subfloat[Deer]{
        \includegraphics[height=0.6in]{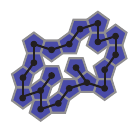}
    }
\caption{Images of holes (at the risk of being psychoanalyzed, the authors couldn't resist the temptation of finding familiar objects in the shape of the holes). Just pentagons with edges that bound holes are shown.}
\label{fig:holes}
\end{figure}

Consequently, up to scale, there are only two distinct types of triangular holes. These correspond to solutions $(a_1,a_2,a_3)=(1,2,2)$ and $(1,1,3)$, resulting in triangles with angles $(36^{\circ}, 72^{\circ}, 72^{\circ})$ and $(36^{\circ}, 36^{\circ}, 108^{\circ})$, respectively. We conjecture that the latter configuration is likely impossible, while the former can only be the one shown in Figure \ref{fig:holes}(a).

The angle types of quadrilateral holes (up to rotation and reflection) can be $(a_1,a_2,a_3, a_4)=$ (1,1,1,7), (1,1,2,6), (1,2,1,6), (1,1,4,4), (1,4,1,4), (1,2,3,4), (1,2,4,3), (1,3,2,4), (2,2,2,4), (2,4,2,4), (2,2,3,3), or (2,3,2,3). We do not know if all of these types are possible or the side length of the holes.

Examples of holes discovered in the simulations are illustrated in Figure \ref{fig:holes}.

\subsection*{Summary}
In the absence of a tessellation, the pentagon cell growth model presents intriguing questions about determining the locus of potential centers and the classification of the types of holes. Our future efforts will focus on resolving the conjectures here presenting and conducting more extensive numerical simulations.

\newpage

\section{Laser Cut Model}

\begin{wrapfigure}[13]{l}{0.35\textwidth}
\centering
\includegraphics[width=0.3\textwidth]{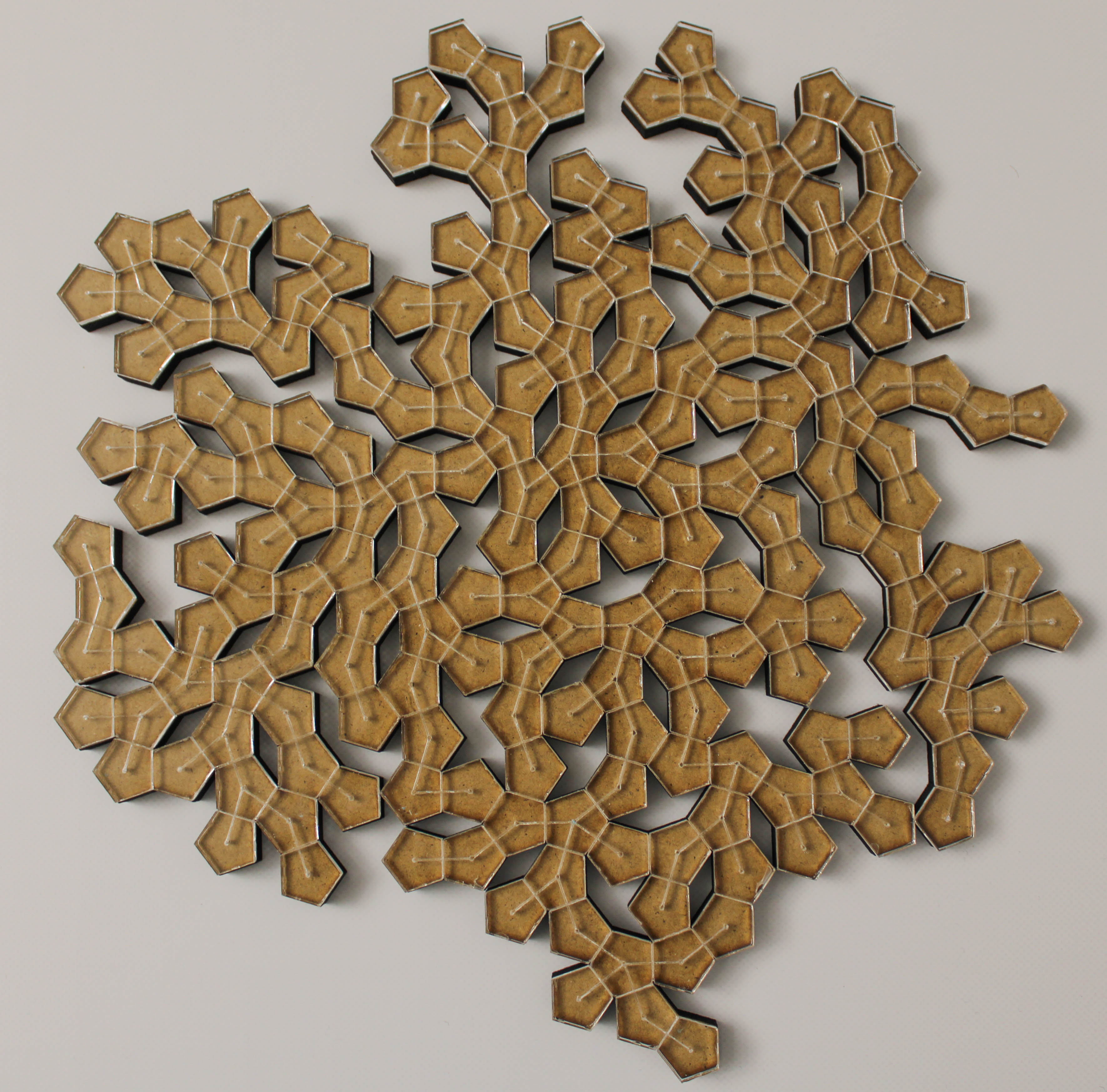}
\caption{\label{a4} The 200 pentagons in wood and acrylic.}
\label{fig:wood}
\end{wrapfigure}
Figure \ref{fig:wood} illustrates both a wooden and acrylic rendition of Figure \ref{pentagons200}(a). We have used the laser cutter Omtech SH-G1390 and the software {\it Adobe Lightroom} and {\it LightBurn}. 
The first tests were carried out exclusively on wood. Currently, the project is being developed using acrylic too. The purpose is to first make the wooden layer, where 
the connected piece is cut (max power: 37\%, velocity: 5 mm/s, number of minutes: 30 approximately); then a second layer on acrylic with the growing tree engraved (max power: 25\%, velocity: 300 mm/s); finally, the connected piece is cut (max power: 25\%, velocity: 10 mm/s, number of minutes in total: 30 approximately). \\
 


{\setlength{\baselineskip}{13pt} 
\raggedright				

\end{document}